\input amstex

\documentstyle{amsppt}
\magnification=1200
\NoBlackBoxes
\def\ls{\vskip.20in}
\def\ss{\vskip.10in}
\def\P{\Bbb P}
\def\bt{$\bullet$}
\def\Q{\Bbb Q}
\def\F{\Cal F}

\centerline{\bf On the Variety of rational Space Curves}
\ss
\centerline{\bf Z. Ran}
\ls

In this paper we study enumeratively the variety $V_{d,n}$ parametrising 
irreducible nonsingular rational curves of degree $d$ in $\P^n, n \geq 3$.  We 
shall give a recursive procedure which computes two sets of numbers associated 
to $V_{d,n}$:  
\roster
\ss
\item"{\bt}"  The Schubert degrees $N_{d,n} (a_1, \ldots, a_k )$, i.e.
the degree of the locus of members of $V_{d,n}$ meeting a generic
collection $A_1, \ldots, A_k$ of linear subspaces of 
respective codimensions $a_a, \ldots, a_k$ in $\P^n$, whenever 
$\sum (a_i - 1) = (n+1) d - (n-3) = \dim V_{d,n}$, i.e. whenever the 
locus in question is finite;
\item"{\bt}"  the linear genera $g_{d,n} (a_1, \ldots, a_k)$, i.e. the 
(geometric) genus of the analogously-defined locus whenever it is 
1-dimensional, i.e. whenever $\sum (a_i - 1) = \dim V_{d,n} - 1$.
\endroster
\ss

Analogous questions for $n=2$ were considered in [1] [2].
As there, the method is completely elementary, involving some geometry
on ruled surfaces.  No quantum-cohomological methods are used, indeed our method
may be viewed as an alternative to those (compare [3] and references therein).  

\subheading{1.  Degrees}
\ss
In what follows we fix $n\geq 3$ and denote by $\bar{V}_d$ the 
closure in the Hilbert scheme of the locus of irreducible nonsingular rational
curves of degree $d$ in $\P^n$.  Let $A_1, \ldots, A_k$ be a 
generic collection of linear subspaces of respective codimensions $a_1,
\ldots, a_k \geq 1$ in $\P^n$.  We denote by 
$$
B_d = B_d (a_{\cdot}) = B_d (A_{\cdot})
$$
the normalization of the locus
$$
\{(C, P_1, \ldots, P_k ) \ : \ C \in \bar{V}_d , P_i  \in C \cap A_i, i =1
, \ldots, k \}
$$
when all $a_i>1$, this
is also the normalization of its projection to $\bar{V}_d$, i.e. the 
locus of degree-$d$ rational curves (and their specializations) meeting
$A_1, \ldots, A_k$; however it will be convenient
to allow some $A_i$ to be hyperplanes.  
Of course if some $a_i>n$ then $N_d(a.)=0$.
Let's call the number of $i$ such that $a_i>1$ the {\it length}
of the condition vector $(a.)$. We have 
$$
\dim B = (n+1) d - (n-3) - \sum (a_i - 1) . 
$$
When this is $0$, we set 
$$
N_d (a_{\cdot}) = \deg B_d (a_{\cdot}) .
$$
Of course $N_d(1,a_2,...)=dN_d(a_2,...),$ so it will suffice to
compute these when all $a_i>1$.

The plan is to get at them via suitable 1-dimensional
$B$'s.  To this end, take a $B = B(A_{\cdot})$ 1-dimensional and let 
$$
\pi : X \to B
$$
be the normalization of the tautological family of rational 
curves, and $f: X \to \P^n$ the natural map.  Then $X$ is smooth and 
each fibre of $\pi$ is either a $\P^1$ or a pair of $\P^1$'s meeting 
transeversely once.  Let $\F$ be the set of components of reducible  
fibres of $\pi$.  Now $X$ comes equipped with a set of 
distinguished sections $s_i = s_{A_i}, i = 1 , \ldots, k$ and  
note that
$$
s_i . s_j = N_d(...,a_i+a_j,...,\hat a_j,...), i \neq j
$$
Also, let $R_i = R_{A_i}$ be the sum of all
fibre components not meeting $s_{A_i}$ and $\F_A \subset \F$ the set of such 
components.  Then $R_{A_i}$ may by blown down, giving rise to a geometric
ruled surface:
$$
b_i : X \to X_i = X_{A_i} = \P (E_{A_i}) ,\tag 1
$$
with sections $\bar{s}_j = b_i (s_j)$ and note that if, say, $i=1$, then 
$$
\bar{s}_2 . \bar{s}_1 = s_2 . s_1 .
$$ 
Now set 
$$
m_i = m_i(a_1,...,a_k) = - s_i^2 , i =1, \ldots, k .
$$
Note that $m_1 = -  \bar{s}_1^2$ too.  But clearly, on a geometric 
ruled surface the
difference of any two sections 
is a sum of fibres hence has self intersection $=0$,
hence 
$$
\align 
0 &= (\bar{s}_1-\bar{s}_2)^2\\
&= s_1^2+s_2^2 + s_2 \cdot R_1 -2s_2\cdot s_1 , \\
&= -m_1 - m_2 + s_2 \cdot R_1   -2s_2\cdot s_1 . 
\endalign 
$$
i.e.
$$
m_1+m_2 = s_2 \cdot R_1 -2N_d(a_1+a_2,a_3,...).  \tag 2 
$$
One consequence of this, already noted and used in [1] is the 

\subheading{2-section lemma 1.1}  If $a_1 = a_2$, then 
$$
s_1^2 = s_2^2 = \frac{-1}{2} s_1 \cdot R_2 + N_d(a_1+a_2,a_3,...)  \tag 3
$$
Indeed if $a_1 = a_2$ then clearly by monodromy $m_1 = m_2$ so (3) 
follows from (2).

For general codimensions we have the 

\subheading{3-section lemma 1.2}  For any 3 distinct distinguished sections 
$s_1, s_2, s_3$ we have 
$$
s_1^2 = \frac{-1}{2} (s_1 \cdot R_2 + s_1 \cdot R_3 - s_2 \cdot R_3 )
+s_1.s_2 + s_1 . s_3 - s_2 . s_3 .  \tag 4 
$$

This follows immediately from (2)
by a suitable linear combination .

Note that by an obvious dimension count the number $k$ of distinguished 
sections on $X$ is always $> 3$, so Lemma 1.2 is always applicable.  Also,
from a recursive standpoint, numbers such as $s_1 \cdot R_2$, having to 
do with reducible curves, are easily computable in terms of $N_d, d' < d$, 
hence may be considered known.  Indeed,
$$
s_1 \cdot R_2 = \sum N_{d_1} (A_{\cdot}^1 , A_1,\P^{s_1} ) 
N_{d_2} (A^2_{\cdot}, 
A_2, \P^{s_2} ). 
\tag 5
$$
the summations being over all $d_1 + d_2 = d, s_1 + s_2 = n $ 
and all decompositions
$A_{\cdot} = (A_1, A_2) \coprod (A_{\cdot}^1) \coprod (A_{\cdot}^2)
$ (as unordered sequences or partitions).  
Each term corresponds to a pair of families of curves
each of degree $d_i$ meeting $(A_1,A_{\cdot}^i)$ and filling up
a locus of codimension $s_i$ and degree $N_{d_i}(A_{\cdot}^i,
A_i,\P^{s_i})$ in $\P^n$, $i=1,2$ ; the two loci meet in a finite
set whose cadinality is given by Bezout's theorem and whose members
correspond with $s_1 \cap R_2$.
\ss
Thus at least when $a_1 + a_2 , a_1 + a_3 , a_2 + a_3 $ are all
$>n$ (which is automatic if $n=3$ but
not otherwise), the 3-section lemma computes $m_1, m_2, m_3$ in terms of
lower-degree data; but even if this condition is not satisfied
(and, say, $a_1, a_2, a_3 >1$), the lemma still computes $m_1$,
say, in terms of data of lower degree {\it or lower length};
we shall use this observation below in constructing a 'length
recursion'.
\ss
Now let $L = f^* \Cal O (1)$ and $F_0$ be a general fibre of $\pi$.  Then we 
have 
$$
L = d s_1 - \sum\limits_{F \in \F_{A_1}} \deg (F) F + x F_0
\tag 6
$$
for some $x \in \Q$ : indeed (6) holds simply because both sides have 
the same value on all fibre components.  To determine $x$, evaluate on $s_1$, 
noting that, by definition,
$$
L\cdot s_1 = N_d (a_1 + 1, a_2, \ldots) .
$$
Thus we have
$$
x = N_d (a_1 + 1, a_2, \ldots) + dm_1 . 
$$
Now let's square (6), noting that, by definition, 
$L^2 = N_s (2, a_1, \ldots)$.  Thus 
$$
N_d (2, a_1, \ldots) = - d^2{m_1} + 2d^2m_1 + 2d N_d (a_1 + 1, a_2, \ldots)
- \sum\limits_{F \in \F_{A_1}} (\deg F)^2
$$
i.e.
$$
N(2,a_1, \ldots) = 2d N_d (a_1 + 1, a_2, \ldots) + d^2 m_1(a.) - 
\sum\limits_{F\in \F_{A_1}} (\deg F)^2 
\tag 7
$$

As above, the sum $\sum\limits_{F\in \F_A} (\deg F)^2$ is easily evaluated in 
terms of $N_{d'} , d'< d$ and may be considered known.  In 
particular, when $a_1 = n$, i.e. $A_1$ is a point,
and moreover $a_2+a_3 > n $, we have  $N_d (a_1 + 1, a_2, \ldots)
= 0$ and that $m_1(a.)$,
via the 3-section lemma, is computable from lower-degree data,
so (7) yields a recursive formula for all the $N_d (2, n, a_2,a_3,
\ldots), a_2 + a_3 > n$, namely
$$
N_d(2,n,a_2,\ldots) = d^2m_1(n,a_2,\ldots) - \sum\limits_{F\in \F_{A_1}}(\deg F)^2
\tag 8
$$
Now take the dot product of (6) with $s_2$, obtaining
$$
N_d(a_1,a_2+1,...) = dN_d(a_1+a_2,...)-\sum\limits_{F
\in \F_{A_1}-\F_{A_2}} (\deg F) + N_d(a_1+1,a_2,...)
+dm_1(a.)
\tag 9
$$
Now to determine $N_d(a.)$ in general we proceed by recursion
on $k$ (as well as $d$), as follows. We may assume all $a_i>1$.
For the smallest possible $k$ (given $d\geq 3$), clearly we may
assume by reordering that
$a_1=a_2=a_3=n$ so from (7) (read backwards) we compute
$N_d(a.)$ from $N_d(2,n-1,n,n,...)=N_d(2,n,n,n-1,...)$
 (and lower-degree data);
but $N_d(2,n,n,n-1,...)$ has already been compted above; this
takes care of the case of smallest $k$. In the general case we
use recursion on $k$. Using (9), we compute $N_d(a.)$ from
$N_d(a_1+1,a_2-1,...)$ and terms of lower degree or length.
Applying (9) repeatedly we compute $N_d(a.)$ in terms of lower-
degree and lower-length data plus $N_d(a_1+a_2-1,1,a_3,...)=
dN_d(a_1+a_2-1,a_3,...)$, itself a lower-length term. This computes
$N_d(a.)$ in general.

\subheading{2.  Genera}
\ss
In this section we fix a sequence $(a_{\cdot})$ giving rise to a smooth (maybe
disconnected) curve $B = B(a_{\cdot})$ and give a formula for the latter's
geometric genus, i.e. for $\deg (K_B)$.  The idea is to consider a `thickening'
$$
B \to B^+ = B (a_{\cdot}^+) = B(A_{\cdot}^+), \dim B^+ = n -1 ,
$$
where $A_{\cdot} = (A_1, \ldots, A_k)$
$$
A_{\cdot}^+ = (A_n^+, \ldots, A_\ell^+ ), \ell \leq k 
$$
$$
\align
A_{\cdot} &\subseteq A_i^+ \subset \P^n\\   
a_i^+&= \dim A_i^+ \leq n-2 . 
\endalign
$$
Consider the diagram (over a neighborhood of the image of 
$B$):
$$
\matrix
X&\longrightarrow&X^+ &\overset{f^+}\to\longrightarrow&\P^n \\
\downarrow&&\downarrow\\
B&\longrightarrow&B^+
\endmatrix
$$
with $f^+$ generically finite.  Let $\rho$ be the ramification
divisor of $f^+$.  Then, as in [2] it is easy to see that
$$
\rho \vert_X = \sum (a_i^+ - 1) s_i . 
$$
On the other hand by Riemann-Hurwitz,
$$
- (n+1) L + \rho = K_{X^+} \vert_X = K_{X/B} + \pi ^*(K_{B^+} \vert_B).
\tag 10 
$$
Now recall the blowing down map $X \to X_1 = X_{A_1} = \P
(E_{A_1})$ (1).  It is easy to see that 
$$
K_{X_1/B} = - 2 s_1 - m_1 F_0 , 
$$
hence
$$
K_{X/B} = - 2 s_1 - m_1 F_0 + R_1
$$
so that by (8),
$$
\pi^* (K_{B^+} \vert_B) = - (n+1) L + \sum (a_i^+ -1) s_i + 2s_1 +
m_1 F_0 - R.
$$
Evaluating on $s_1$, we conclude
$$
\deg (K_{B^+} \vert_B) = -(n+1) N_d (a_1 + 1, a_2, \ldots) - 
a_1^+ m_1. 
$$
By the adjunction formula, 
$$
\deg K_B = \deg (K_{B^+} \vert_B) + \deg N_{B/B^+}
$$
On the other hand, if we set
$$
\align
B_i &= B(a_1, \ldots, a_i^+, \ldots, a_k) \quad i = 1, \ldots, \ell\\
&= B (a_1, \ldots, \hat{a}_i, \ldots)\quad i > \ell
\endalign
$$
then clearly
$$
\deg N_{B/B^+} = \sum\limits_{i=1}^k \deg N_{B/B_i} , 
\tag 11 
$$
so it suffices to evaluate $\deg N_{B/B_i}$.

\subheading{Case 1}  $i \leq \ell$.

We then have a Cartesian diagram
$$
\matrix
s_{A_i} & \longrightarrow & A_i \\
\cap&&\cap\\
s_{A_i^+} &\longrightarrow & A_i^+
\endmatrix
$$
from which clearly
$$
\pi^* N_{B/B_i} = N_{s_{A_i}/s_{A_i^+}} = f^* ((a_i - a_i^+) \Cal O(1)) 
$$
hence
$$
\deg N_{B/B_i} = (a_i - a_i^+) N(a_1, \ldots, a_i + 1, \ldots ) \tag 12
$$

\subheading{Case 2}  $i > \ell$

We then have a Cartesian diagram
$$
\matrix
&&s_{A_i} &\longrightarrow& A_i \\
&&\cap \\
B&\longleftarrow & X && \cap \\
\downarrow&&\downarrow \\
B_i&\longleftarrow& X_i &\longrightarrow& \P^n
\endmatrix
$$
from which as above 
$$
\deg N_{s_{A_i}/X-i} = a_i N_d (a_1, \ldots, a_i + 1, \ldots ) 
.
$$
On the other hand
$$
\align
\deg N_{s_{A_i}/X_i} & = \deg N_{s_{A_i}/X} + \deg N_{X/X_i}\vert _{s_{A_i}}\\
&= - m_i + \deg \pi^* N_{B/B_i} .
\endalign
$$
Consequently
$$
\deg N_{B/B_i} = a_i N_d (a_1, \ldots, a_i +1, \ldots) + m_i . \tag 13
$$
Putting (9) (10) (11) together, we have computed 
$\deg (K_B)$, as claimed.
\subheading{References}

1. Z. Ran : 'Bend, break and count' (alg-geom 9704004)

2. Z. Ran: 'Bend, break and count II: elliptics, cuspidals, linear genera'
(alg-geom 9708013) 

3. W. Fulton and R. Pandharipande: 'Notes on stable maps and quantum cohomology'
Mittag-Leffler Report No. 4, 1996/97.
\enddocument